\documentclass[
final
]{dmtcs-episciences}

\usepackage[utf8]{inputenc}
\usepackage{subfigure}
\usepackage{amsmath}
\usepackage{amsthm}
\newtheorem{theorem}{Theorem}[section]
\newtheorem{Lemme}{Lemma}[section]
\theoremstyle{Remarque}
\newtheorem{Remarque}{Remark}[section]
\newtheorem{Proposition}{Proposition}[section]
\newtheorem{Corollaire}{Corollary}[section]
\theoremstyle{Exemple}
\newtheorem{Exemple}{Example}[section]
\theoremstyle{definition}
\newtheorem{definition}{Definition}[section]

\usepackage[round]{natbib}

\author[Moussa Barro et al.]{Moussa Barro\affiliationmark{1}%\thanks{I am fully supported.}
  \and K. Ernest Bognini\affiliationmark{2}%\thanks{And they are, too!}
  \and Boucar{\'e} Kient{\'e}ga\affiliationmark{3}}
\title[On the Study of cellular automata on modulo-recurrent words ]{On the study of cellular automata on modulo-recurrent words}
\affiliation{
Nazi BONI University, Bobo-Dioulasso, Burkina Faso\\
Joseph KI-ZERBO University, Ouagadougou, Burkina Faso\\
Daniel Ouezzin COULIBALY University, Dedougou, Burkina Faso}
\begin{document}
% This is only used if you are compiling for a volume before vol 25
% \publicationdetails{VOL}{2015}{ISS}{NUM}{SUBM}
% This is the new form of collecting the data, starting with vol 25
\publicationdata{vol. 27:3}{2025}{29}{10.46298/dmtcs.10380}{2022-11-28; 2022-11-28; 2025-01-03; 2025-04-22; 2025-08-25; 2025-12-09}{2025-12-09}
\maketitle

\begin{abstract}
\vspace{0.3cm}
In this paper, we study a class of cellular automata (CA) called stable cellular automata (SCA) that preserve stability by reflection, modulo-recurrent, and richness. After applying these automata to Sturmian words, we determine some of their combinatorial properties. Next, we calculate the classical and palindromic complexity functions of these words. Finally, we demonstrate that these words are $2$-balanced and establish their abelian complexity function.
\end{abstract}

%DMTCS is an open access scientific is implemented by the
%\emph{episcience} platform, see \cite{berthaud:hal-01002815} for an
%overview of the strategy. It combines high scientific and editorial
%quality with an open access policy. It is priceless, neither authors
%nor readers pay money for the access. Access is granted by giving
%episcience an irrevocable license to publish the articles, the
%copyright remains with the authors. The platform itself is run by
%French government services that do their best to warrant continuous
%access and a high quality of service.
%
%This document describes the use of the \texttt{dmtcs-episcience.cls}
%document class. It has to be used \emph{for all DMTCS publications}.

\section{Introduction}
\label{sec:in}
A cellular automaton is a series of cells evolving using a precise set of rules, resulting in a new generation of cells. These automata were introduced in \cite{b0} with the objective to realize some dynamic systems capable to model complex sefl-reproduction phenomena. Later, in the 1970s, the concept was popularized by the work of John Horton Conway with his famous \emph{game of life} on two-dimensional cellular automata (CA). Thus, CA have become a multidisciplinary field of study going from physics to computer science and from biology to mathematics. 

 Modulo-recurrent words are recurrent words in which any factor appears in all positions modulo all integers. For instance, we have Sturmian words and Champernowne word.  These words were introduced in \cite{KT-f}  and intensively sudied in \cite{BKT, CKT, b23}.     
         
A palindrome is a word which is the same when we read left to right or from right to left. The study of palindromic factors in combinatorics on words enables us to characterize certain infinite words (see \cite{b7,b9,b24}). 
                                       
Given a finite or infinite word, the complexity function $\textrm{c}$ calculates the number of distinct factors of a given length in the latter. Studying this function allows us to characterize certain families of infinite words \cite{b6}. This concept has also been used to establish various characterizations and properties of Sturmian words (see \cite{b27, b5, b20,b28,b12,b16}). Depending on the specific factors included in a word and whether they are finite or infinite, several types of complexity functions can be distinguished, such as palindromic and abelian functions. The palindromic complexity function computes the number of distinct palindromic factors of a given length within a word. The abelian complexity function counts the number of Parikh vectors for each length of the word. This function was intensively studied in \cite{BKT-p, BKT}. These two notions allow us to characterize Sturmian words \cite{b11}. 
In this work, we study the combinatorial properties of infinite words obtained through the application of cellular automata (CA). This paper is organized as follows: First, we provide some definitions and notation. Then, in Section 2, we review properties of Sturmian and modulo-recurrent words. In Section 3, we apply CA to infinite words and demonstrate that these automata preserve properties such as modulo-recurrence and periodicity. Section 4 defines a class of CA called stable cellular automata (SCA), establishing that they preserve stability by reflection and richness. Lastly, in Section 5, we discuss the combinatorial study of words obtained by applying these SCA to Sturmian words.
%\clearpage
\section{Preliminaries}
\label{sec:first}
\subsection{Definitions and notation}
\label{subsec:first}
An alphabet $\mathcal{A}$, is a finite set whose elements are called letters. A word is a finite or infinite sequence of elements over $\mathcal{A}$. We denote by $\mathcal{A}^\ast $, the set of finite words over $\mathcal{A}$ and $\varepsilon$ the empty word. For all $u\in \mathcal{A}^*$, $|u|$ denotes the length of $u$ and $|u|_x$ for all $x$ in $\mathcal{A}$, the number of occurrence of $x$ in $u$. A word $u$ of length $n$ constitued by a single letter $x$ is simply denoted $u=x^n$; by convention $x^0=\varepsilon $. Let $u=u_1u_2 \dotsb u_n$ be a finite word with  $u_i\in \mathcal{A}$ for all $i\in\left\{1,2,\dots ,n\right\}$. The word $\overline{u}=u_n\dotsb$ $u_2u_1$ is called the reflection of $u$. Given two finite words $u$ and $v$ then, we have $\overline{uv}=\overline{v}$
 $\overline{u}$. The word $u$ is called palindrome if $\overline{u}=u$.

    We denote by $\mathcal{A}^{\omega}$ (respectively, $ \mathcal{A}^\infty=\mathcal{A}^*\cup \mathcal{A}^{\omega}$ ), the set of infinite (respectively, finite and infinite) words.
    
     An infinite word $\textbf{u}$ is ultimately periodic if there are two words $v\in \mathcal{A}^*$ and $w\in \mathcal{A}^+$ such that $\textbf{u}=vw^\omega$. The word $\textbf{u}$ is called recurrent if each of its factors appears an infinite number of times. Moreover, $\textbf{u}$ is called periodic if $v=\varepsilon$. The $n$-th power for some finite word $w$ is denoted by $w^n$.
  
  Let $\textbf{u}\in \mathcal{A}^\infty$ and $v\in \mathcal{A}^*$. We say that $v$ is a factor of $\textbf{u}$ if there exists $u_1\in \mathcal{A}^*$ and $\textbf{u}_2\in \mathcal{A}^\infty$ such that $\textbf{u}=u_1v\textbf{u}_2$. In other words, we say that $\textbf{u}$ contains $v$.
We also say that $u_1$ is a prefix of $\textbf{u}$ and we note $u_1=\text{Pref}_{|u_1|}(\textbf{u})$. In particular, if $\textbf{u} \in \mathcal{A}^*$ then $\textbf{u}_2$ is called suffix of $\textbf{u}$. 
    
  Let $w$ be a factor of an infinite word $\textbf{u}$ and $x$ a letter in $\mathcal{A}$. Then, $L_n(\textbf{u})$ denotes the set of factors of length $n$ of $\textbf{u}$ and $L(\textbf{u})$ that of all factors of $\textbf{u}$. The letter $x$ is called a left (respectively, right) extension of $w$ if $xw$ (respectively, $wx$) belongs to $L(\textbf{u})$. Let us denote by $\partial^-w$ (respectively, $\partial^+w$) the number of left (respectively, right) extension of $w$ in $\textbf{u}$. When $\partial^+w=k$ with $k\geq 2$, $w$ is called a right $k$-prolongeable factor. In the same way, we can define the notion of left $k$-prolongeable factors. A factor $w$ of $\textbf{u}$ is called a right (respectively, left) special if $\partial^+w>1$ (respectively, $\partial^-w>1$). Any factor that is both right and left special is called a bispecial factor.
  
  Given an infinite word $\textbf{u}$. Let $\texttt{N}=\{1,2,\cdots \}$ denote  the set of natural numbers and let $\mathbb{N}=\{0,1,2,\cdots\}$ denote the set of whole numbers. The map of $\mathbb{N}$ into $\texttt{N}$ defined by  $\text{c}_\textbf{u}(n) = \# L_n(\textbf{u})$, where $\# L_n(\textbf{u})$ denotes the cardinality of $L_n(\textbf{u})$, is called the complexity function of $\textbf{u}$. This function is related to the special factors by the following relation (see \cite{b8} for details):
    $$\hspace{0 cm} \text{c}_\textbf{u}(n+1)-\text{c}_\textbf{u}(n)= \displaystyle\sum_{w\in L_n(\textbf{u})} (\partial^+(w)-1). $$

    We denote the set of palindromic factors of length $n$ by $\text{Pal}_n(\textbf{u})$ and the set of all palindromic factors of $\textbf{u}$ by $\text{Pal}(\textbf{u})$. The palindromic complexity function of $\textbf{u}$, denoted $\text{c}^{\text{pal}}_\textbf{u}$, is the map from $\mathbb{N}$ to $\texttt{N}$ that counts the number of distinct palindromic factors of length $n$ contained in $\textbf{u}$:
$$\text{c}^{\text{pal}}_\textbf{u}(n) = \# \left\{w\in L_n(\textbf{u}) : \overline{w}=w \right\}.$$ 
If, for all $w\in L(\textbf{u})$, we have $\overline{w}\in L(\textbf{u})$, then $\textbf{u}$ is called stable by reflection.

  Let $w$ be a factor of the infinite word $\textbf{u}$ over the alphabet $\mathcal{A}_q=\{a_{1}, a_{2}, \cdots, a_{q}\}$. Then, the $q$-uplet
     $\chi (w)=(| w |_{a_{1}}, | w|_{a_{2}}, \cdots, | w |_{a_{q}})$ is called the Parikh vector of $w$. The set of Parikh vectors of factors of length $n$ in $\textbf{u}$ is denoted by:  $$\chi_{n}(\textbf{u})=\{\chi(w): w\in L_{n}(\textbf{u})\}.$$ 
     The abelian complexity function of $\textbf{u}$ is defined as the map from $\mathbb{N}$ to $\texttt{N}$ by: $$\text{c}^{\text{ab}}_\textbf{u}(n)= \# \chi_{n}(\textbf{u}).$$

  The window complexity function of $\textbf{u}$ is the map,  $\text{c}^\text{win}_\textbf{u}$ from $\mathbb{N}$ to $\texttt{N}$, defined by $$\text{c}^\text{win}_\textbf{u}(n)=\#\left\{u_{kn}u_{kn+1}\cdots u_{n(k+1)-1} : k\geq 0\right\}.$$ 

The shift $S$, is an application defined over $\mathcal{A}^\omega$  which erases the first letter of a given word. For instance, if $\textbf{u}=u_0u_1u_2u_3\cdots$ then $S(\textbf{u})=u_1u_2u_3\cdots$. A substitution $\varphi$ is a map of $\mathcal{A}^*$ into itself such that $\varphi (uv)=\varphi(u)\varphi(v)$, for 
any $u$, $v\in \mathcal{A}^*$.

\subsection{Sturmian words and modulo-recurrent words}
\label{subsec:Two}
In this subsection, we will review some properties of Sturmian and modulo-recurrent words that will be used later. We will consider the alphabet to be the set of two letters, $\mathcal{A}_{2}=\left\{a,b\right\}$, because Sturmian words are binary.
\begin{definition}\label{D1}
  An infinite word $\textbf{u}$ over $\mathcal{A}_{2}$  is called Sturmian if for any whole number $n,\ \text{c}_\textbf{u}(n)=n+1$.
  \end{definition}
  The most well-known Sturmian word is the  probably famous Fibonacci word. It is the fixed point of the substitution $\varphi$ defined over $\mathcal{A}_{2}^\ast$ by:
   $$\varphi(a)=ab \ \textrm{and} \ \varphi(b)=a.$$
   It is noted:
  $$\mathbf{F}=\displaystyle \lim_{n\rightarrow \infty}\varphi^n(a).$$

 \begin{definition}
   A Sturmian word is called $a$-Sturmian (respectively, $b$-Sturmian) when it contains $ a^2 $ (respectively, $ b^2$).
\end{definition}  
 
\begin{definition}
  An infinite word $\textbf{u}$ over $\mathcal{A}_{2}$  is called quasi-Sturmian if there exist $k\in \texttt{N}$ and $n_0\in \mathbb{N}$ such that for any $n\geq n_0,\ \text{c}_\textbf{u}(n)=n+k$.
  \end{definition}
  
  \begin{definition}
  Let $\textbf{u}$ be an infinite word over $\mathcal{A}$ and $x\in \mathcal{A}$. The set of return words for a letter $x$ in $\textbf{u}$ is defined by: $$ret_{\textbf{u}}(x)=\left\{xw\in L(\textbf{u}): |w|_x=0, xwx\in L(\textbf{u})\right\}.$$
  \end{definition}

\begin{definition} A word $\textbf{u}=u_0u_1u_2 \dotsb$ is called modulo-recurrent if any factor of $\textbf{u}$ appears in all positions modulo $i$ for all $i\geq 1$.
\end{definition}
\begin{definition} Let $w$ be a factor of an infinite word $\textbf{u}$. We say that $w$ is a window factor if $w$ appears in $\textbf{u}$ at a position which is a multiple of $|w|$.
\end{definition}

\begin{Proposition}\cite{CKT}\label{propo-comp-mod-rec}
Let $\textbf{\emph{u}}$ be a modulo-recurrent word. Then, for all integers $n$, the set of window factors of length $n$ in $\textbf{\emph{u}}$ is equal to $L_n(\textbf{\emph{u}})$.
\end{Proposition}
       
 \begin{definition} An infinite word $\textbf{u}$ is called $\alpha$-balanced if $\alpha$ is the smallest integer such that for any pair ($v$, $w$) of factors in $\textbf{u}$ of the same length and for every letter $x$ in $\mathcal{A}$, we have: 
    \begin{center}
      $||v|_x-|w|_x|\leq \alpha$.
      \end{center} 
      
      If $\alpha =1$, then $\textbf{u}$ is called simply balanced.
      
\begin{Lemme}\label{L}\cite{10}
Let $\textbf{\emph{u}}$ be an infinite unbalanced word over the binary alphabet $\mathcal{A}$. Then, there exists a factor, $v_{1}$ of $L(\textbf{\emph{u}})$ and two distinct letters, $x$ and $y$, such that $xv_{1}x$ and $yv_{1}y$ are in  $L(\textbf{\emph{u}})$.
\end{Lemme}

     \end{definition}
     The following theorem provides a classic characterization of Sturmian words.
     
     \begin{theorem}\label{theo-stur} \cite{b11, b9} Let $\textbf{\emph{u}}$ be an infinite binary word. The following assertions are equivalent: 
     \begin{enumerate}
     \item $\textbf{\emph{u}}$ is Sturmian.
     
     \item $\textbf{\emph{u}}$ is non-ultimately periodic and balanced.

      \item For all $n\geq 1, \hspace{0.3cm} \text{c}^{\emph{ab}}_\textbf{\emph{u}}(n)=2$.
     
     \item For all $n\geq 0, \hspace{0.3cm}$ $$ \text{c}^{\emph{pal}}_\textbf{\emph{u}}(n) = \left \{
 \begin{array}{l}
 \hspace{0cm}1 \hspace{0.3cm} \text{if}\ n\ \text{is even} \hspace{0.1cm} \\
 \hspace{0cm} 2 \hspace{0.33cm} \text{otherwise}. 
 \end{array}
 \right. $$
 
     \end{enumerate}

     \end{theorem}$\vspace{0.5cm}$
     
     \begin{theorem}\label{stur-puis}\cite{b16} 
     Let $\textbf{\emph{v}}$ be an $a$-Sturmian word over $\mathcal{A}_2$. Then, there exists a Sturmian sequence $(\epsilon_i)_{i\geq 1}$ over the alphabet $\left\{0,1\right\}$ and integer $l$  such that $\textbf{\emph{v}}$ is written:    
     $$\textbf{\emph{v}}=a^{l_0}ba^{l+\epsilon_1}ba^{l+\epsilon_2}ba^{l+\epsilon_3}b \dotsb,$$  where $l_o\leq l+1$.
     \end{theorem}
 It is proved in \cite{b23} that Sturmian words are modulo-recurrent.
 
\begin{theorem}\label{stur-mod}\cite{CKT}
Let $\textbf{\emph{u}}$ be an infinite recurrent word. The following assertions are equivalent:
\begin{enumerate}
\item $\textbf{\emph{u}}$ is a modulo-recurrent word. 
\item For all integers $n\geq 1,\  \text{c}^{\emph{win}}_\textbf{\emph{u}}(n)= \text{c}_\textbf{\emph{u}}(n)$.
\end{enumerate}
  
\end{theorem}
\section{Cellular automata (CA)}
\label{sec:CA}
We suppose in the following that you write a paper since yIn this section, we define a class of CA that we apply to infinite words. For the remainder of this paper, $\textbf{u}\in \mathcal{A}^{\infty}$ and $F$ is a CA defined over $\mathcal{A}^r$.
  \begin{definition}
Let $\mathcal{A}$ and $\mathcal{B}$ be two alphabets. Then, a map $F:\mathcal{A}^*\longrightarrow\mathcal{B}^*$ is a cellular automaton (CA) if there exists an integer $r\geq 1$ and a morphism $f :\mathcal{A}^r \longrightarrow  \mathcal{B}$ such that:
 $$ \left \{
 \begin{array}{l}
 \hspace{0cm}F(w)=\varepsilon \hspace{3.2cm} \text{if} \hspace{0.3cm} |w|< r \hspace{0.1cm} \\
 \hspace{0cm}F(xyz)= f(xy)F(yz) \hspace{0.95cm} \text{if} \hspace{0.3cm} x\in \mathcal{A}, \ y\in \mathcal{A}^{r-1},\  z\in \mathcal{A}^{\infty}.
 \end{array}
 \right. $$  
  \end{definition}

  From this definition, we have the following remark.

  \begin{Remarque} \label{rem-long}  
 For all finite words, we have that $|F(w)|=|w|-r+1$ when $|w|\geq r$.
  \end{Remarque}
    In the following, we will say that $F$ is injective if $f$ is.

  \begin{Proposition}\label{fac-conserv}  
 Let $\textbf{\emph{u}},\textbf{\emph{v}}\in \mathcal{A}^{\infty}$ and $F$ be a CA. Then, we have:
 \begin{enumerate}
\item If $u_1\in L(\textbf{\emph{u}})$ then $F(u_1)\in L(F(\textbf{\emph{u}}))$.
\item If $w\in L(F(\textbf{\emph{u}}))$ then there exists $u_1\in L(\textbf{\emph{u}})$ such that $w=F(u_1)$.
 
\item If $ L(\textbf{\emph{u}})\subset L(\textbf{\emph{v}})$ then $L(F(\textbf{\emph{u}}))\subset L(F(\textbf{\emph{v}}))$.
 
 \item $L(F(\textbf{\emph{u}}))= F(L(\textbf{\emph{u}}))$.

 \end{enumerate}
 
\end{Proposition}
\begin{proof}
 \begin{enumerate}
\item This follows from the definition of $F$.
\item Suppose that $w\in L_n(F(\textbf{u}))$, then we can write $w=y_1y_2\cdots y_n$. Therefore, there exists $u_1\in L(\textbf{u})$ with $u_1=x_1x_2\cdots x_{n+r-1}$ such that $y_1=F(x_1x_2\cdots x_{r})$, $y_2=F(x_2x_3\cdots x_{r+1})$, $\cdots$ and $y_n=F(x_nx_{n+1}\cdots x_{n+r-1})$ because $F$ is a function. Consequently, $w=F(u_1)$.
\item Suppose that $ L(\textbf{u})\subset L(\textbf{v})$. Let $w\in  L(F(\textbf{u}))$, then by 2., there exists $u_1\in L(\textbf{u})$ such that $w=F(u_1)$. Since $ L(\textbf{u})\subset L(\textbf{v})$ and  $u_1\in L(\textbf{u})$, then $u_1\in L(\textbf{v})$. Therefore, by 1., we have that $F(u_1)\in L(F(\textbf{v}))$, i.e, $w\in L(F(\textbf{v}))$. Therefore, $L(F(\textbf{u}))\subset L(F(\textbf{v}))$.
 
 \item  If $w\in  L(F(\textbf{u}))$, then there exists $u_1\in L(\textbf{u})$ such that $w=F(u_1)$. Since $ F(u_1)\in F(L(\textbf{u}))$, we have that $w\in F(L(\textbf{u}))$. Therefore, $L(F(\textbf{u}))\subset F(L(\textbf{u}))$.
Conversely, if $w\in  F(L(\textbf{u}))$, then there exists $u_1\in L(\textbf{u})$ such that $w=F(u_1)$. Since $u_1\in L(\textbf{u})$, then we have $F(u_1)\in L(F(\textbf{u}))$ by 1, so $w\in L(F(\textbf{u}))$. Consequently, $F(L(\textbf{u}))\subset L(F(\textbf{u}))$. Therefore, $L(F(\textbf{u}))= F(L(\textbf{u}))$.  
 \end{enumerate}
\end{proof}

 \begin{theorem}\label{cor-inj} For any infinite word $\textbf{\emph{u}}$, we have: $$\emph{c}_{F(\textbf{\emph{u}})}(n)\leq \emph{c}_\textbf{\emph{u}}(n+r-1), \forall\ n\in \mathbb{N}.$$
   Moreover, $\emph{c}_{F(\textbf{\emph{u}})}(n)= \emph{c}_\textbf{\emph{u}}(n+r-1)$, for any $ n\in \texttt{N} $, if $F$ is injective.
  \end{theorem}
  \begin{proof}\\
$\bullet$ According to Remark \ref{rem-long}, any factor of $F(\textbf{u})$ of length $n$ comes from a factor of length $n+r-1$ of $\textbf{u}$. Moreover, we have $\# L_{n}(F(\textbf{u}))\leq \# L_{n+r-1}(\textbf{u})$. Therefore, $\text{c}_{F(\textbf{u})}(n)\leq \text{c}_\textbf{u}(n+r-1)$, for all $n\in \mathbb{N}$.\\  
$\bullet$ Suppose that $F$ is injective. In that case, each factor of $F(\textbf{u})$ of length $n$ comes from only one factor of length $n+r-1$ of $\textbf{u}$. Thus, we obtain $\# L_{n}(F(\textbf{u}))\geq \# L_{n+r-1}(\textbf{u})$. Hence, $\text{c}_{F(\textbf{u})}(n)= \text{c}_\textbf{u}(n+r-1)$, for all $n\in \texttt{N}$. \end{proof}

  \begin{theorem}\label{mod}
  Let $\textbf{\emph{u}}$ be an infinite word and $F$ be a CA.  
  \begin{enumerate}
   \item If $\textbf{\emph{u}}$ is modulo-recurrent, then $F(\textbf{\emph{u}})$ is modulo-recurrent. 
  
  \item If $F$ is injective and $F(\textbf{\emph{u}})$ is modulo-recurrent, then $\textbf{\emph{u}}$ is modulo-recurrent.
  \end{enumerate}
 
  \end{theorem}
  \begin{proof}
  \begin{enumerate}
  \item Suppose that $\textbf{u}$ is modulo-recurrent. Let $w\in L(F(\textbf{u}))$. By Proposition \ref{fac-conserv}, there exists $u_1\in L(\textbf{u})$ such that $w=F(u_1)$. Additionally, if $u_1$ appears at position $j$ in $\textbf{u}$, then $F(u_1)$ appears at position $j$ in $F(\textbf{u})$. Since $\textbf{u}$ is modulo-recurrent, then $u_1$ appears in every position modulo $ i$ in $\textbf{u}$, for all $i\geq 1$. Therefore, $F(u_1)$ also appears in every position modulo $i$ in $F(\textbf{u})$, for all $i\geq 1$. Thus, $F(\textbf{u})$ is also modulo-recurrent.
  \item Sippose that $F$ is injective and that $F(\textbf{u})$ is modulo-recurrent. If $u_1\in L(\textbf{u})$, then Proposition \ref{fac-conserv} tells us that $F(u_1)\in L(F(\textbf{u}))$. Since $F$ is injective, the factors $u_1$ and $F(u_1)$ appear at the same positions in both $\textbf{u}$ and $F(\textbf{u})$. Since $F(\textbf{u})$ is modulo-recurrent, the factor $F(u_1)$ appears in every position modulo $i$ in $F(\textbf{u})$, $\forall i\geq 1$. Therefore $u_1$ also appears in every position modulo $i$ in $\textbf{u}$, for all $i\geq 1$. Consequently, $\textbf{u}$ is also modulo-recurrent.
  \end{enumerate}
 \end{proof}  
  
  \begin{Lemme} Let $\textbf{\emph{u}}$ be an infinite word. The following assertions then hold.
  \begin{enumerate}
  \item If $\textbf{\emph{u}}$ is periodic then $F(\textbf{\emph{u}})$ is periodic.
  \item If $F$ is injective and $F(\textbf{\emph{u}})$ periodic then $\textbf{\emph{u}}$ is periodic.   
  \end{enumerate}
  \end{Lemme}
  \begin{proof}
    \begin{enumerate}
  \item Suppose that $\textbf{u}$ is periodic. Then, there exists a finite word $u_1$ such that $\textbf{u}=u_1^\omega.$ As a result,
  $$\textbf{u}=u_1^\omega\Longrightarrow  F(\textbf{u})=v_{1}^\omega,$$
  where  $v_1= F(\text{Pref}_{|u_1|+r-1}( \textbf{u}))=\text{Pref}_{|u_1|}( F(\textbf{u}))$. Therefore, $F(\textbf{u})$ is periodic.
  \item Suppose that $F$ is injective and that $F(\textbf{u})$ is periodic. Then, there exists a factor $v_{1}$ of $F(\textbf{u})$ such that $F(\textbf{u})=v_{1}^{\omega}$. Since $F$ is injective, there is a unique factor $w\in L_{\vert v_{1}\vert+r-1}( \textbf{u})$ such that $F(w)=v_1$. By setting $u_1=\text{Pref}_{\vert v_1\vert}(w)$, we obtain that $\textbf{u}=u_1^{\omega}$. Therefore, $\textbf{u}$ is periodic.
  \end{enumerate}
\end{proof}
\begin{Remarque}
Let $F$ be an injective CA. If $\textbf{\emph{u}}$ and $F(\textbf{\emph{u}})$ are two periodic words, then they have the same minimal period.
\end{Remarque}
   \begin{theorem}
  Let $\textbf{\emph{u}}$ be a recurrent word and $F$ a CA such that:
$$ \left \{
 \begin{array}{l}
 \hspace{0cm}F(x_{1}y_{1})= F(x_{2}y_{2})\hspace{1cm} if\ x_{1}= x_{2} \\
 \hspace{0cm}F(x_{1}y_{1})\neq F(x_{2}y_{2}) \hspace{1cm} \text{otherwise},
 \end{array}
 \right.$$
where $x_{1}, x_{2}\in \mathcal{A}$ and $\ y_{1}, y_{2}\in \mathcal{A}^{r-1}$. Then, $F(\textbf{\emph{u}})$ is balanced if and only if $\textbf{\emph{u}}$ is balanced.
  \end{theorem}
  \begin{proof}
 Consider two alphabets $\mathcal{A}$ and $\mathcal{B}$, that each have exactly two letters. 

Suppose that $\textbf{u}$ is balanced and $F(\textbf{u})$ is unbalanced. Since $F(\textbf{u})$ is unbalanced, there exists a factor $v_1\in L(F(\textbf{u}))$ and $x'_{1},x'_{2}\in\mathcal{B}$ and $x'_1\neq x'_2$ such that $x'_1v_1x'_1$, $x'_2v_1x'_2\in L(F(\textbf{u}))$, by Lemma \ref{L}. According to Proposition \ref{fac-conserv}, there are two factors $u_1,\ u_2\in L(\textbf{u})$ such that $x'_1v_1x'_1=F(u_1)$ and $x'_2v_1x'_2=F(u_2)$. Thus, we can write $u_1=au_1'a\delta_1$ and $u_2=bu_1'b\delta_2$ where $\delta_1,\ \delta_2\in L_{r-1}(\textbf{u})$ and $a,b\in\mathcal{A}$ with $a\neq b$. Consequently, $au_1'a,\ bu_1'b\in L(\textbf{u})$. We get a contradiction because $\textbf{u}$ is balanced.
  
  Conversely, suppose that $F(\textbf{u})$ is balanced and that $\textbf{u}$ is unbalanced. Since $\textbf{u}$ is unbalanced, there exists a factor $u_1\in L(\textbf{u})$ and $a,\ b\in \mathcal{A}$ with $a\neq b$ such that $au_1a,\ bu_1b\in L_n(\textbf{u})$. Additionally, there are $\delta_1,\ \delta_2\in L_{r-1}(\textbf{u})$ such that $au_1a\delta_1,\ bu_1b\delta_2\in L_{n+r-1}(\textbf{u})$. Furthermore, $F(au_1a\delta_1)=x'_1v_1x'_1$ and $F(bu_1b\delta_2)=x'_2v_1x'_2$ are some factors of $L(F(\textbf{u}))$. This contradicts our hypothesis. 
  
  We obtain the desired equivalence from all of the above.
  \end{proof}

  \begin{theorem}\label{spec} Let $\textbf{\emph{u}}$ be a recurrent word and $F$ a CA over $\mathcal{A}$ such that: 
$$ \left \{
 \begin{array}{l}
 \hspace{0cm}F(x_{1}y_{1})= F(x_{2}y_{2})\hspace{1cm} if\ x_{1}= x_{2} \\
 \hspace{0cm}F(x_{1}y_{1})\neq F(x_{2}y_{2}) \hspace{1cm} \text{otherwise},
 \end{array}
 \right.$$
where $x_{1}, x_{2}\in \mathcal{A}$ and $ y_{1}, y_{2}\in \mathcal{A}^{r-1}$. Let $u_1\in L_n(\textbf{\emph{u}})$ with $n\geq r$. If $F(u_1)$ is left special factor of $F(\textbf{\emph{u}})$ then there is a left special factor of $\textbf{\emph{u}}$ (called $u_{2}$) such that $F(u_2)=F(u_1)$.
\end{theorem}

\begin{proof} 
Suppose that $u_1\in L_n(\textbf{u})$ with $n\geq r$ such that $v_1=F(u_1)$ is left special factor of $F(\textbf{u})$. Then, $x'_1v_1$ and $x'_2v_1$ are in $L(F(\textbf{u}))$ where $x'_{1},x'_{2}\in\mathcal{B}$ and $x'_1\neq x'_2$. According to Proposition \ref{fac-conserv}, there exists $u_2\in L(\textbf{u})$ and $a,b\in\mathcal{A}$; $a\neq b$ such that $x'_1v_1=F(au_2)$ and $x'_2v_1=F(bu_2)$. Therefore, $au_2, bu_2 \in L(\textbf{u})$. As a result, $u_2$ is left special factor in $\textbf{u}$ and $F(u_2)=F(u_1)$.  
\end{proof}

 \begin{theorem}\label{spec} Let $\textbf{\emph{u}}$ be a recurrent word and $F$ a CA over $\mathcal{A}$ such that: 
$$ \left \{
 \begin{array}{l}
 \hspace{0cm}F(y_{1}x_{1})= F(y_{2}x_{2})\hspace{1cm} if\ x_{1}= x_{2} \\
 \hspace{0cm}F(y_{1}x_{1})\neq F(y_{2}x_{2}) \hspace{1cm} \text{otherwise},
 \end{array}
 \right.$$
where $x_{1}, x_{2}\in \mathcal{A}$ and $ y_{1}, y_{2}\in \mathcal{A}^{r-1}$. Let $u_1\in L_n(\textbf{\emph{u}})$ with $n\geq r$. If $F(u_1)$ is right special factor of $F(\textbf{\emph{u}})$ then there is a right special factor of $\textbf{\emph{u}}$ (called $u_{2}$)  such that $F(u_2)=F(u_1)$.
\end{theorem}

\begin{proof} 
Suppose that $u_1\in L_n(\textbf{u})$ with $n\geq r$ such that $v_1=F(u_1)$ is right special factor of $F(\textbf{u})$. Then, $v_1x'_1$ and $v_1x'_2$ are in $L(F(\textbf{u}))$ where $x'_{1},x'_{2}\in\mathcal{B}$ and $x'_1\neq x'_2$. According to Proposition \ref{fac-conserv}, there exists $u_2\in L(\textbf{u})$ and $a,b\in\mathcal{A}$; $a\neq b$ such that $v_1x'_1=F(u_2a)$ and $v_1x'_2=F(u_2b)$. Therefore, $u_2a, u_2b \in L(\textbf{u})$. As a result, $u_2$ is right special factor in $\textbf{u}$ and $F(u_2)=F(u_1)$. 
\end{proof}

\section{Stable cellular automata (SCA)}
\label{sec:SCA}
In this section, we will study a class of cellular automata called stable cellular automata (SCA).

\begin{definition}
A cellular automaton $F$ is invariant if, for any infinite word $\textbf{u}$,  $F(\textbf{u})=\textbf{u}$.
\end{definition}

The following proposition provides a characterization  of an invariant cellular automaton.

\begin{Proposition}\label{eq1} Let $F$ be a CA. Then, $F$ is invariant if and only if $F(xy)=x$ for all $x\in \mathcal{A}\ \text{and} \ y\in \mathcal{A}^{r-1}$.

\end{Proposition}
\begin{proof} Let $\textbf{u}$ be an infinite word in the form $\textbf{u}=x_0x_1x_2\cdots$ with $x_i\in \mathcal{A}$, for all $i\in \mathbb{N}$. Then, 
$$F(\textbf{u})=F(x_0x_1\cdots x_{r-1})F(x_1x_2\cdots x_{r})F(x_2x_3\cdots x_{r+1})\cdots.$$ As a result, the following equivalences hold:
\begin{align*}
F\ \text{is invariant}&\Longleftrightarrow F(\textbf{u})=\textbf{u},\hspace{0.2cm} \forall \textbf{u}\in \mathcal{A}^{\omega}\\
&\Longleftrightarrow F(x_0x_1\cdots x_{r-1})F(x_1x_2\cdots x_{r})F(x_2x_3\cdots x_{r+1})\cdots=x_0x_1x_2x_3 \cdots\\
&\Longleftrightarrow F(x_ix_{i+1}\cdots x_{i+r-1})=x_i, \hspace{0.2cm}\forall x_i\in \mathcal{A},\  i\in \mathbb{N}\\
&\Longleftrightarrow F(xy)=x,\ \forall\ x\in \mathcal{A},\ y\in \mathcal{A}^{r-1}.
\end{align*}
\end{proof}
  
\begin{Lemme}\label{Lem-Ech}
Let $F$ be an invariant CA over $\mathcal{A}^r_2$ and $E$ be the exchange map defined over $\mathcal{A}_2$. Then, $F \circ E=E \circ F$.
\end{Lemme}
\begin{proof} Let $\textbf{u}$ be an infinite word over $\mathcal{A}_2$ and $u_1\in L(\textbf{u}) $. Then, we distinguish two cases.
\begin{enumerate}[{Case }1:]
\item $|u_1|< r$. Then, we have $F(u_1)=\varepsilon$. As a result, $E (F(u_1))=\varepsilon$. Additionally, $|E(u_1)|<r$. Thus, $F(E(u_1))=\varepsilon$. Hence, $F(E(u_1))=\varepsilon=E (F(u_1))$.
\item $|u_1|\geq r$. Let $w$ be a factor of length $r$ in $u_1$. Then, we can write in the form $w=xw_1$, with $x\in\mathcal{A}_2$. We have $E(w)=E(x)E(w_1)$. Thus, by Proposition \ref{eq1}, we have $F(E(w))=E(x)$. Furthemore, $F(w)=x$, i.e, $E(F(w))=E(x)$. As a result, we obtain $F(E(w))=E(F(w))$. It follows that $F(E(u_1))=E(F(u_1))$.
\end{enumerate}
In all cases, $F \circ E=E \circ F$.
\end{proof}
\begin{definition} Let $F$ be a CA. Then, $F$ is called stable if $F(\overline{w})=F(w)$, for all $w\in \mathcal{A}^r$.
\end{definition}

\begin{Lemme}\label{Lem-stable}
Let $\textbf{\emph{u}}$ be an infinite word and $F$ a SCA. Then, for all $u_1\in L(\textbf{\emph{u}})$, we have  $F(\overline{u_1})=\overline{F(u_1)}$.
\end{Lemme}
\begin{proof}
 Let $u_1\in L(\textbf{u})$. Then, we distinguish two cases.
 \begin{enumerate}[{Case }1:]
\item $|u_1|< r$. Then, $F(u_1)=F(\overline{u_1})=\varepsilon=\overline{F(u_1)}$.

\item $|u_1|\geq r$. Then, we have $u_1\in L_{n+r}(\textbf{u})$ and $F(u_1)=F(w_0)\cdots F(w_n)$ where $w_i=\text{Pref}_r(S^i(u_1))$, for all $i\in \{0,\dots ,n\}$. 

Additionally, we have:
\begin{align*}
\overline{F(u_1)}&=\overline{F(w_0)F(w_1)\cdots F(w_n)}\\
&=\overline{F(w_n)}\cdots \overline{F(w_1)}\ \overline{F(w_0)}\\
&=F(w_n)\cdots F(w_1)F(w_0),\ \text{because}\ F(w_i)\in \mathcal{A}, \ \text{for all}\ i\in \{0\dots,n\}.
\end{align*}
Furthermore,  we have $F(\overline{ u}_1)=F(\overline{ w}_n)\cdots F(\overline{ w}_1)\ F(\overline{ w}_0)$. Since $F$ is stable we have, for all $i\in \{0,\dots,n\}$,  $F(w_i)=F(\overline{ w}_i)$ with $w_i\in L_r(\textbf{u})$. Thus, $F(\overline{u}_1)=F(w_n)\cdots F(w_1)F(w_0)$.\\
Hence, $\overline{F(u_1)}=F(\overline{u}_1)$.
\end{enumerate}
\end{proof}
\begin{theorem}\label{Theo-stable}  
Let $\textbf{\emph{u}}$ be an infinite word and $F$ a SCA. 
\begin{enumerate}
\item If $\textbf{\emph{u}}$ is stable by reflection, then $F(\textbf{\emph{u}})$ is stable by reflection.
\item If $F$ is injective and $F(\textbf{\emph{u}})$ is stable by reflection, then $\textbf{\emph{u}}$ is stable by reflection.
\end{enumerate}
\end{theorem}
\begin{proof} 
\begin{enumerate}
\item Suppose that $\textbf{u}$ is stable by reflection. Let $u_1\in L_{n+r}(\textbf{u})$. Then, by Proposition \ref{fac-conserv} we have $F(u_1)\in L_{n+1}(F(\textbf{u}))$. By writing $F(u_1)=F(w_0)F(w_1)\cdots F(w_n)$ with $w_i=\text{Pref}_r(S^i(u_1)$, for all $i\in \{0,\dots,n\}$. Then, By Lemma \ref{Lem-stable},  $\overline{F(u_1)}=F(\overline{u}_1)$. However, $F(\overline{u}_1)\in L(F(\textbf{u}))$. Thus, $\overline{F(u_1)}\in L(F(\textbf{u}))$. Hence, $F(\textbf{u})$ is stable by reflection.

\item Assume that $F$ is injective and $F(\textbf{u})$ is stable by reflection. Let $v_1\in L_{n+1}(F(\textbf{u}))$ then, there exists $u_1\in L_{n+r}(\textbf{u})$ such that $v_1=F(u_1)$. As a result, $F(u_1)\in L_{n+1}(F(\textbf{u}))$ and $\overline{F(u_1)}\in L_{n+1}(F(\textbf{u}))$ because $F$ is stable. However, by Lemma \ref{Lem-stable},  $\overline{F(u_1)}=F(\overline{u}_1)$, i.e, $F(\overline{u}_1) \in L(F(\textbf{u}))$. Since $F(\overline{u}_1) \in L(F(\textbf{u}))$ and $F$ is injective, then we have $\overline{u}_1\in L(\textbf{u})$. Hence, $\textbf{u}$ is stable by reflection.
\end{enumerate}
\end{proof}

  \begin{Corollaire}\label{pal1} Let $\textbf{\emph{u}}$ be an infinite word that is stable under reflection and let $F$ be an injective SCA. Then, a factor $u_1$ of $\textbf{\emph{u}}$ is a palindromic factor if and only if $F(u_1)$ is a palindromic factor.
  \end{Corollaire}
  \begin{proof}
Suppose that $u_1$ is a palindromic factor of $\textbf{u}$. Then, by Lemma \ref{Lem-stable}, we have $F(\overline{u}_1)=\overline{F(u_1)}$. Additionally, $F(\overline{u}_{1})=F(u_1)$. Hence, $\overline{F(u_1)}=F(u_1)$.
  
  Reciprocaly suppose that $\overline{F(u_1)}=F(u_1)$. As $F(\overline{u}_1)=\overline{F(u_1)}$, by Lemma \ref{Lem-stable} we get  
$F(\overline{u}_1)=F(u_1)$. Since $F$ is injective, we deduce that $\overline{u_1}=u_1$. 
\end{proof}

 \begin{Corollaire}\label{pal2} Let $\textbf{\emph{u}}$ be an infinite word stable by reflection and $F$ an injective SCA. Then, for all $n\in \texttt{N}$, we have:
  $$\emph{c}^{\emph{pal}}_{F(\textbf{\emph{u}})}(n)=\emph{c}^{\emph{pal}}_{\textbf{\emph{u}}}(n+r-1).$$
  \end{Corollaire}
\begin{proof}
 Use Theorem \ref{cor-inj} and Corollary \ref{pal1}.
\end{proof}

  \begin{definition}
  Let $\textbf{u}\in \mathcal{A}^{\infty}$. We say that $\textbf{u}$ is rich if any factor $w$ of $\textbf{u}$, has exactly $|w|+1$ distinct palindromic factors including the empty word. 
 \end{definition}

\par In \cite{bucci1}, the authors make the following useful remark on rich infinite words. Indeed, it shows necessary and sufficient condition of an infinite word to be rich.  
  \begin{Remarque}\label{rq-rich}
A word $\textbf{\emph{u}}$ is rich if and only if any factor of $\textbf{\emph{u}}$ is also rich. 
  \end{Remarque} 
  
 The result below in \cite{bucci} characterizes the rich words stable by reflection.
 \begin{theorem}\label{Bucci}
Let $\textbf{\emph{u}}$ be an infinite word such that the set of its factors is stable under reflection. Then, $\textbf{\emph{u}}$ is rich if and only if, for all $n \in \mathbb{N}$, we have:
  $$\emph{c}^{\emph{pal}}_\textbf{\emph{u}}(n)+\emph{c}^{\emph{pal}}_\textbf{\emph{u}}(n+1)=\emph{c}_{\textbf{\emph{u}}}(n+1)-\emph{c}_\textbf{\emph{u}}(n)+2.$$
 \end{theorem} 

\begin{Corollaire}\label{C}
Sturmian words are rich.
\end{Corollaire}
\begin{proof}\\
 We obtained the result by first applying Definition \ref{D1} and Theorem \ref{theo-stur}, and then Theorem \ref{Bucci}. 
\end{proof}
  
   \begin{theorem}\label{riche}
Let $\textbf{\emph{u}}$ be an infinite word that is stable under reflection and let $F$ be an injective SCA. If $\textbf{\emph{u}}$ is rich, then $F(\textbf{\emph{u}})$ is also rich. 
  \end{theorem}
 \begin{proof}
Since $F$ is injective and $\textbf{u}$ is stable under reflection, respectively by Theorem \ref{cor-inj} and Corollary \ref{pal2}, for all $n\in \texttt{N}$, $\text{c}_{F(\textbf{u})}(n)=\text{c}_{\textbf{u}}(n+r-1)$ and $\text{c}^{pal}_{F(\textbf{u})}(n)=\text{c}^{pal}_{\textbf{u}}(n+r-1)$.\\
Additionally,  $\text{c}_{F(\textbf{u})}(n+1)-\text{c}_{F(\textbf{u})}(n)=\text{c}_\textbf{u}(n+r)-\text{c}_\textbf{u}(n+r-1)$, for all $n\in \texttt{N}$. Moreover, since $\textbf{u}$ is recurrent and stable by reflection, then by Theorem \ref{Bucci}, we have:
  \begin{align*}
  \textbf{u}\ \text{rich} \Longrightarrow \text{c}_\textbf{u}(n+r)-\text{c}_\textbf{u}(n+r-1)+2 &=\text{c}_\textbf{u}^{\text{pal}}(n+r)+\text{c}^{\text{pal}}_\textbf{u}(n+r-1), \forall n\in \texttt{N}\\
&=\text{c}^{\text{pal}}_{F(\textbf{u})}(n+1)+\text{c}^{\text{pal}}_{F(\textbf{u})}(n).
%\forall n\in \texttt{N}\ \text{by %injectivity of}\ F.
  \end{align*}
Consequently, $\text{c}_{F(\textbf{u})}(n+1)-\text{c}_{F(\textbf{u})}(n)+2=\text{c}^{\text{pal}}_{F(\textbf{u})}(n+1)+\text{c}^{\text{pal}}_{F(\textbf{u})}(n)$, $\forall n\in \texttt{N}$.\\
 Thus, $F(\textbf{u})$ is rich. 
\end{proof}
\section{Cellular automata and Sturmian words}
\label{sec:CA and Sturmian words}

In this section, we apply the CA to Sturmian words. 
  \begin{definition} A CA, $F$ is called Sturmian CA if, for any Sturmian word $\textbf{u}$, the word $F(\textbf{u})$ is also Sturmian. 
  \end{definition}
   
\begin{Exemple}   
Let us consider $\textbf{\emph{u}}$ a Sturmian word over $\mathcal{A}_2$. Then, for the CA defined by:  
$$ \left \{
 \begin{array}{l}
 \text{for all}\ x\in \mathcal{A}_{2}\ \text{and}\ y\in \mathcal{A}_{2}^{r-1} \\
H(xy)=x\\ 
G(xy)=E(x),
 \end{array}
 \right.$$ 
we have $H(\textbf{\emph{u}})=\textbf{\emph{u}}$ and $G(\textbf{\emph{u}})=E(H(\textbf{\emph{u}}))=H(E(\textbf{\emph{u}}))=E(\textbf{\emph{u}})$ which are respectively fixed point and Sturmian word where $E$ is the exchange map defined over $\mathcal{A}_2$.
    \end{Exemple} 
Note that $H$ and $G$ are SCA.
   
  \subsection{Classical and window complexity}
  Let $\textbf{v}$ be a Sturmian word over $\mathcal{A}_2$. Then, $\textbf{v}$ is can be written as   $\textbf{v}=a^{l_0}ba^{l+\epsilon_1}ba^{l+\epsilon_2}ba^{l+\epsilon_3}b \dotsb$ with $(\epsilon_i)_{i\geq 1}$ a Sturmian sequence over $\left\{0,1\right\}$ and $l_o\leq l+1$, by Theorem \ref{stur-puis}.  \\

    Consider the SCA $F$ defined over $\mathcal{A}_2^{l+1}$ by:
   
     $$F(w)= \left \{
 \begin{array}{l}
 \hspace{0cm}a \hspace{0.2cm} if \hspace{0.2cm} w=a^{l+1} \\
 \hspace{0cm} b \hspace{0.2cm} otherwise.
 \end{array}
 \right. $$
 
 Then, we get:
    
    $$F(\textbf{v})=\left \{
 \begin{array}{l}
 ab^{l+1}a^{\epsilon_1}b^{l+1} a^{\epsilon_2}b^{l+1}a^{\epsilon_3}b^{l+1}a^{\epsilon_4}b^{l+1} \dotsb \hspace{0.45cm} if \hspace{0.2cm} l_o=l+1 \\ 
  b^{l_0+1}a^{\epsilon_1}b^{l+1} a^{\epsilon_2}b^{l+1}a^{\epsilon_3}b^{l+1}a^{\epsilon_4}b^{l+1} \dotsb  \hspace{0.5cm} otherwise.
 \end{array}
 \right. $$
 
 Clearly, $F$ is a SCA.
 
 Let $k_0$ be the maximum power of $a^{l}b$ (respectively $a^{l+1}b$) in $\textbf{v}$ if the sequence $(\epsilon_i)_{i\geq 1}$ is $0$-Sturmian (respectively $1$-Sturmian). The number $k_0$ is called the maximal power of $a^{l}b$ in $\textbf{v}$ if $(a^{l}b)^{k_0}\in L(\textbf{v})$ and $(a^{l}b)^{k_0+1}\notin L(\textbf{v})$.    
     
     Let $n_0=k_0(l+1)$ in the following.

 \begin{Lemme} \label{ret}  
 The sets of return words for the letters of $F(\textbf{\emph{v}})$ are given by: $$ret_{F(\textbf{\emph{v}})}(a)=\left\{ab^{n_0-l-1}, \ ab^{n_0} \right\}\hspace{1cm} \text{and}\hspace{1cm} ret_{F(\textbf{\emph{v}})}(b)=\left\{b, \ ba \right\}.$$ 
\end{Lemme}
\begin{proof} Note that $F(\textbf{v})$ can be written as follows: $$F(\textbf{v})=\left \{
 \begin{array}{l}
 ab^{n_0-\epsilon'_1}ab^{n_0-\epsilon'_2}ab^{n_0-\epsilon'_3}ab^{n_0-\epsilon'_4}a \dotsb \hspace{1.6cm} if \hspace{0.2cm} l_o=l+1 \\ 
 b^{l_0+1}a^{\epsilon_1}b^{n_0-\epsilon'_1}ab^{n_0-\epsilon'_2}ab^{n_0-\epsilon'_3}ab^{n_0-\epsilon'_4}a \dotsb  \hspace{0.5cm} otherwise
 \end{array}
 \right. $$ where $(\epsilon'_i)_{i\geq 1}$ is a sequence over $\left\{0,l+1\right\}$. Thus, the return words of the letter $a$ are $\left\{ab^{n_0-l-1}, \ ab^{n_0} \right\}$. Similary, the return words of the letter $b$ are $\left\{b, \ ba \right\}.$ 
 \end{proof}

We study the combinatorial properties, such as the classical and window complexities, of the word $F(\textbf{v})$.
  \begin{Proposition}  
  For all $n\geq 0,\ \emph{c}^{\emph{win}}_{F(\textbf{\emph{v}})}(n)=\emph{c}_{F(\textbf{\emph{v}})}(n).$
  \end{Proposition}
  \begin{proof}
 Since $\textbf{v}$ is modulo-recurrent, $F(\textbf{v})$ is also modulo-recurrent, by Theorem \ref{mod}. Therefore, by Theorem \ref{stur-mod}, $\text{c}^{\text{win}}_{F(\textbf{v})}(n)=\text{c}_{F(\textbf{v})}(n)$ . 
 \end{proof}  
\begin{Lemme}\label{Lem1-auto-stur}
Let $v_1$ be a factor of $F(\textbf{\emph{v}})$ such that $|v_1|>n_0$. 
Then, there exists exactly one $w \in L(\textbf{\emph{v}})$ such that $F(w) = v_1$.
\end{Lemme}  

\begin{proof}
  Note that the letter $a$ has a single preimage, $a^{l+1}$ in $L(\textbf{v})$ by $F$. Additionally, any factor of $F(\textbf{v})$ of length  strictly greater than $n_0$ contains at least one occurrence of $a$. Therefore, this factor must come from a single preimage in $L(\textbf{v})$. 
  \end{proof} 
 \begin{Remarque}\label{rq-unique}
 If $w_1,w_2\in L_{n+r-1}(\textbf{\emph{v}})$ with $w_1\neq w_2$ and $n>n_0$, then $F(w_1)\neq F(w_2)$.
  \end{Remarque}
 
  \begin{theorem}\label{cc}
The classical complexity function of the word $F(\textbf{\emph{v}})$ is given by: $$\emph{c}_{F(\textbf{\emph{v}})}(n)= \left \{
 \begin{array}{l}
 \hspace{0cm}n+1 \hspace{2.5cm} if \ n\leq n_0-l \\
 \hspace{0cm}2n-n_0+l+1 \hspace{0.72cm} if \ n_0-l< n\leq n_0 \\
 \hspace{0cm} n+l+1 \hspace{1.9cm} if \ n>n_0. 
 \end{array}
 \right. $$ 
  \end{theorem}
\begin{proof}
Let us proceed by disjunction of cases according to the length $n$ of the factors.
\begin{enumerate}[{Case }1:]
\item $1\leq n\leq n_0-l$. Observe that the factors of $F(\textbf{v})$ of length $n$ contain at most one occurrence of the letter $a$. Then, we have $L_n(F(\textbf{v}))=\left\{b^n,\ b^iab^{n-i-1} :\ i=0,\dots, n-1 \right\}.$ Therefore, $ \text{c}_{F(\textbf{v})}(n)=n+1$, for all $n\leq n_0-l$. Let us observe that for all $ n\leq n_0-l-1$, $b^n$ is the only right special factor of length $n$ of $F(\textbf{v})$. Similarly, the right special factor of $n_0-l$ in $F(\textbf{v})$ are $b^{n_0-l}$ and $ab^{n_0-l-1}$. Hence, we obtain the following equalities: 
  \begin{align*}
  \text{c}_{F(\textbf{v})}(n_0-l+1)&=\text{c}_{F(\textbf{v})}(n_0-l)+2\\
  &=(n_0-l+1)+2\\
  &=n_0-l+3.
  \end{align*}
\item  $n_0-l+1\leq n\leq n_0$. Note that the factors of $F(\textbf{v})$ of length $n$ contain at most two occurrences of the letter $a$. Then, we have:\\  
$L_n(F(\textbf{v}))=\{b^n,\ b^iab^{n-i-1}, \ b^jab^{n_0-l-1}ab^{n-n_0+l-1-j} :\ i=0,\dots, n-1; \\
\hspace*{4cm} 0\leq j \leq n-n_0+l-1 \}.$ \\
  As a result, we get:
  \begin{align*}
\text{c}_{F(\textbf{v})}(n)&=1+n+(n-n_0+l-1+1)\\
 &=2n-n_0+l+1.
  \end{align*}
\item $n> n_0$. Then, any factor of $F(\textbf{v})$ of length $n$ comes from only one factor of length $n+r-1$ in $\textbf{v}$, by Lemma \ref{Lem1-auto-stur}. By applying Theorem \ref{cor-inj}, we obtain the following equalities:
  \begin{align*}
 \text{c}_{F(\textbf{v})}(n)&=\text{c}_\textbf{v}(n+r-1)\\
  &=n+r\\
  &=n+l+1.
  \end{align*}
  \end{enumerate}
  \end{proof}
  \begin{Remarque}
The word $F(\textbf{\emph{v}})$ is a quasi-Sturmian word.
 
  \end{Remarque}

\subsection{Palindromic properties}
In this subsection, we will study the palindromic complexity function and the palindromic richness of $F(\textbf{v})$.

\begin{Lemme}\label{Lem2-auto-stur} Let $\emph{v}$ be Sturmian and $F$ be a SCA. Then, any palindromic factor of $F(\textbf{\emph{v}})$ of length greater than $n_0$ comes from a palindromic factor of $\textbf{\emph{v}}$.
\end{Lemme}  
\begin{proof}
Let $v_1$ be a palindromic factor of $F(\textbf{v})$ such that $|v_1|> n_0$. Then, by Lemma \ref{Lem1-auto-stur}, $v_1$ comes from only one factor $u_1$ of $\textbf{v}$. Thus, since |$v_1|> n_0$, we can apply the same reasoning as in Theorem \ref{Theo-stable} and Corollary \ref{pal1} to show that $u_1$ is a palindromic factor.
\end{proof}
 
 \begin{theorem}\label{cp}
The palindromic complexity function of a word $F(\textbf{\emph{v}})$ is given by:
 \begin{enumerate}
\item If $n \leq n_0-l$,
 $$\emph{c}^{\emph{pal}}_{F(\textbf{\emph{v}})}(n)= \left \{
 \begin{array}{l}
 \hspace{0cm}1 \hspace{0.2cm} if \hspace{0.2cm} n \hspace{0.1cm} is \hspace{0.1cm} even  \\
 \hspace{0cm}2 \hspace{0.2cm} otherwise. 
 \end{array}
 \right. $$ 
 
\item If $n_0-l< n \leq n_0$, 
\begin{itemize}
\item for $n_0$ even, we have: $ \emph{c}^{\emph{pal}}_{F(\textbf{\emph{v}})}(n)= \left \{
 \begin{array}{l}
 \hspace{0cm}1 \hspace{0.2cm} if \hspace{0.2cm} l \hspace{0.1cm} and \hspace{0.1cm} n \hspace{0.1cm} are \hspace{0.1cm} even \\
 \hspace{0cm}2 \hspace{0.2cm} if \hspace{0.2cm} l \hspace{0.1cm} is \hspace{0.1cm}  odd\\
 \hspace{0cm}3 \hspace{0.2cm} otherwise, 
 \end{array}
 \right. $ 
\item for $n_0$ odd, we have: $\emph{c}^{\emph{pal}}_{F(\textbf{\emph{v}})}(n)=2.$
\end{itemize}
\item If $n> n_0$,
 
 $$\emph{c}^{\emph{pal}}_{F(\textbf{\emph{v}})}(n)= \left \{
 \begin{array}{l}
 \hspace{0cm}1 \hspace{0.2cm} if \hspace{0.2cm} n+l \hspace{0.1cm} is \hspace{0.1cm} even \\
 \hspace{0cm}2 \hspace{0.2cm} otherwise.
 \end{array}
 \right. $$ 

 \end{enumerate}
  
 \end{theorem}
\begin{proof}
 \begin{enumerate}
\item If $ n \leq n_0-l$, then we have: $$L_n(F(\textbf{v}))=\left\{b^n,\ b^iab^{n-i-1} : i=0,1,\dots, n-1 \right\}.$$
Therefore, $b^iab^{n-i-1}$ is a palindromic factor of length $n$ if and only if $i=n-i-1$, i.e, $n=2i+1$. As a result, we obtain:  $$\text{Pal}_n(F(\textbf{v}))= \left \{
 \begin{array}{l} 
  \left\{b^n \right\}
  \hspace{3cm}if\hspace{0.1cm} n \hspace{0.1cm}is \hspace{0.1cm} even\\
 \left\{b^n; \hspace{0.1cm}b^{\frac{n-1}{2}} ab^{\frac{n-1}{2}} \right\} \hspace{0.9cm} otherwise.
 \end{array}
 \right. $$
 
\item If $n_0-l< n\leq n_0$ then:
 $$ \hspace*{-1.2cm}L_n(F(\textbf{v}))=\left\{b^n,\ b^iab^{n-i-1}, \ b^jab^{n_0-l-1}ab^{n-n_0+l-1-j} : i=0,1,\dots, n-1; \  0\leq j \leq n-n_0+l-1 \right\}.$$
Therefore, $b^jab^{l+1}ab^{n-n_0+l-1-j}$ is a palindromic factor of length $n$ of $F(\textbf{v})$ if and only if $j=n-n_0+l-1-j$. It follows that, $n+l=2j+n_0+1$. Thus, let us now reason according to parity of $n_0$ to ensure that the word $b^jab^{l+1}ab^{n-n_0+l-1-j}$ be a palindromic factor of $F(\textbf{v})$: 
\begin{itemize}
\item for $n_0$ even, as $n+l=2j+n_0+1$, we deduce that $l$ and $n$ are different parities. As a conseqence, we obtain:  $$ \hspace*{-1.2cm} \text{Pal}_n(F(\textbf{v}))= \left \{
 \begin{array}{l} 
  \left\{b^n \right\} \hspace{7.7cm} if\hspace{0.1cm} l \hspace{0.1cm}and\hspace{0.1cm} n\hspace{0.1cm} are\hspace{0.1cm} even\\
 \left\{b^n;\hspace{0.1cm} b^{\frac{n-1}{2}} ab^{\frac{n-1}{2}} \right\} \hspace{5.6cm}if\hspace{0.1cm}l \hspace{0.1cm} and\hspace{0.1cm} n\hspace{0.1cm} are\hspace{0.1cm} odd\\
 \left\{b^n;\hspace{0.1cm} b^{\frac{n-n_0+l-1}{2}}ab^{n_0-l-1}ab^{\frac{n-n_0+l-1}{2}} \right\} \hspace{2.5cm}if\hspace{0.1cm}l \hspace{0.1cm}is \hspace{0.1cm}odd\hspace{0.1cm} and\hspace{0.1cm} n \hspace{0.1cm}is \hspace{0.1cm}even\\
 \left\{b^n;\hspace{0.1cm} b^{\frac{n-1}{2}} ab^{\frac{n-1}{2}};\hspace{0.1cm} b^{\frac{n-n_0+l-1}{2}}ab^{n_0-l-1}ab^{\frac{n-n_0+l-1}{2}} \right\} \hspace{0.5cm}if\hspace{0.1cm}l \hspace{0.1cm}is \hspace{0.1cm}even\hspace{0.1cm} and\hspace{0.1cm} n \hspace{0.1cm}is \hspace{0.1cm}odd.
 \end{array}
 \right. $$
\item for $n_0$ odd, since $n+l=2j+n_0+1$ then the integers $n$ and $l$ have  different parities. Thus:  $$\text{Pal}_n(F(\textbf{v}))= \left \{
 \begin{array}{l} 
 \left\{b^n;\hspace{0.1cm} b^{\frac{n-1}{2}} ab^{\frac{n-1}{2}} \right\} \hspace{5cm}if\hspace{0.1cm} n\hspace{0.1cm} is\hspace{0.1cm} odd\\
 \left\{b^n;\hspace{0.1cm} b^{\frac{n-n_0+l-1}{2}}ab^{n_0-l-1}ab^{\frac{n-n_0+l-1}{2}} \right\} \hspace{1.9cm}if\hspace{0.1cm} n \hspace{0.1cm}is \hspace{0.1cm}even.
\end{array}
 \right. $$
\end{itemize} 
\item If $n> n_0$. Then firstly, it is known that any factor of length $n$ of $F(\textbf{v})$ comes from a factor of length $n+l$ of $\textbf{v}$. Secondly, by Lemma \ref{Lem2-auto-stur}, any palindromic factor of $F(\textbf{v})$ of length $n> n_0$, comes from only one palindromic factor of $\textbf{v}$. Additionally, by applying the Theorem \ref{theo-stur}, we deduce that $F(\textbf{v})$ has a palindromic factor of length $n$ if $n+l$ is even and two otherwise.
 \end{enumerate}
 \end{proof}
 
 \begin{Corollaire}
The word $F(\textbf{\emph{v}})$ is rich.
 \end{Corollaire}
\begin{proof}
 This follows from Corollary \ref{C} and Theorem \ref{riche}.
\end{proof} 

\subsection{Abelian complexity function} 
 In this subsection, we will determine the balance, Parikh vectors and abelian complexity function of $F(\textbf{v})$.
 
\begin{Proposition}\cite{BKT-p}\label{pab-binaire}
Let $\textbf{\emph{u}}$ be an infinite $\beta$-balanced word over $\left\{a,b\right\}$. Then, for all integers $n$, we have: $$\emph{c}^{\emph{ab}}_\textbf{\emph{u}}(n)\leq \beta+1.$$
\end{Proposition}

\begin{Lemme}\label{l-eq} 
Let $u_1, u_2\in L_n(\textbf{\emph{v}})$. If $|u_1|=|u_2|$, then we have $||u_1|_{a^{l+1}}-|u_2|_{a^{l+1}}|\leq 2$.
\end{Lemme}
\begin{proof}  
 Let $u_1, u_2\in L_n(\textbf{v})$. Suppose that there exists a minimal integer $n$ such that \\
$|u_1|_{a^{l+1}}=|u_2|_{a^{l+1}}+3$. Then, we can write $u_{1}$ and $u_{2}$ in the form:
 $$u_1=a^{l+1}u'_1a^{l+1}\ \text{and}\ u_2=\alpha_1u'_2\alpha_2,\ \text{with}\ \alpha_1, \alpha_2\in \left\{a^{i}ba^{l-i}:\ i=0,1,\dots,l \right\}.$$ Consequently, we have:
  $$|u'_1|=|u'_2|,\ |u'_1|_{a^{l+1}}=|u'_2|_{a^{l+1}}+1,\ |u'_1|_b=|u_1|_b\ \text{and}\ |u_2|_b=|u'_2|_b+2.$$
   We obtain  $||u_1|_{b}-|u_2|_{b}|= ||u'_1|_{b}-(|u'_2|_{b}+2)|\leq 1$ becauce $\textbf{v}$ is balanced. Thus,\\ $|u'_1|_{b}-(|u'_2|_{b}+2)\in\left\{-1,0,1  \right\}$. We distinguish the following cases:

$\bullet$ if $|u'_1|_{b}-(|u'_2|_{b}+2)=-1$ then, $|u'_1|_{b}-|u'_2|_{b}=-3 $.

$\bullet$ if $|u'_1|_{b}-(|u'_2|_{b}+2)=1$ then, $|u'_1|_{b}-|u'_2|_{b}=-1 $.

$\bullet$ if $|u'_1|_{b}-(|u'_2|_{b}+2)=0$ then, $|u'_1|_{b}-|u'_2|_{b}=-2 $.
 
The three cases are impossible. Therefore, we deduce that $||u_1|_{a^{l+1}}-|u_2|_{a^{l+1}}|\leq 2$.
\end{proof}

 \begin{Lemme}\label{eq} 
 The word $F(\textbf{\emph{v}})$ is $2$-balanced.
\end{Lemme}
\begin{proof}
 Let $v_1, v_2\in L_n(F(\textbf{v}))$. We distinguish the following cases according to the length $n$ of factors $v_1$ and $v_2$.
\begin{enumerate}[{Case }1:]
\item $1\leq n \leq n_0-l$. Then, since
 $L_n(F(\textbf{v}))=\left\{b^n,\ b^iab^{n-i-1} :\ i=0,1,\dots, n-1 \right\}$.

Hence, we have $||v_1|_x-|v_2|_x|\leq 1$, $\forall x\in \left\{a,b\right\}$.

\item $n_0-l+1\leq n\leq n_0$. Then, we have:
 $$L_n(F(\textbf{v}))=\{b^n,\ b^iab^{n-i-1},\ b^jab^{n_0-l-1}ab^{n-n_0+l-1-j} : i=0,\dots, n-1; \ 0\leq j \leq n-n_0+l-1\}.$$ 
Therefore, $||v_1|_x-|v_2|_x|\leq 2$, $\forall x\in \left\{a,b\right\}$.\\
 
\item $n\geq n_0+1$. Since the word $\textbf{v}$ is Sturmian, we note that for all factors $u_1$ and $u_2$ of $\textbf{v}$, we have: 
 
 $\hspace{2cm}|u_1|=|u_2| \ \Longrightarrow\  ||u_1|_{a^{l+1}}-|u_2|_{a^{l+1}}|\leq 2$ by Lemma \ref{l-eq}.\\ 
 Additionally, $F(a^{l+1})=a$ then, $|u_1|_{a^{l+1}}=|F(u_1)|_a$. Since $v_1, v_2\in L_n(F(\textbf{v}))$ then, there are two factors $u_1, u_2\in L_{n+l}(\textbf{v})$ such that $v_1=F(u_1)$ and $v_2=F(u_2)$. Since $||u_1|_{a^{l+1}}-|u_2|_{a^{l+1}}|\leq 2$, then we have $||F(u_1)|_a-|F(u_2)|_a|\leq 2$, i.e, $||v_1|_a-|v_2|_a|\leq 2$. As $|v_1|_a+|v_1|_b=|v_1|$, then we obtain $||v_1|_b-|v_2|_b|\leq 2$. Therefore, $F(\textbf{v})$ is $2$-balanced. 
\end{enumerate}
\end{proof}

 \begin{theorem}\label{ca}
 The abelian complexity function of $F(\textbf{\emph{v}})$ is given by the following formulas for all $n\in \texttt{N}$ by:
 
 \begin{enumerate}
\item For $1\leq n \leq n_0-l$, $\emph{c}^{\emph{ab}}_{F(\textbf{\emph{v}})}(n)=2$.
\item For $n_0-l+1\leq n\leq n_0$,  $\emph{c}^{\emph{ab}}_{F(\textbf{\emph{v}})}(n)=3$.
\item For $n> n_0 $,  $\text{c}^{ab}_{F(\textbf{\emph{v}})}(n)\in \left\{2,3 \right\}$.
 \end{enumerate}
  
\end{theorem}  
\begin{proof}
We distinguish the following cases based on the length $n$ of the factors.
\begin{enumerate}
\item If $1\leq n \leq n_0-l$. Then, since
 $L_n(F(\textbf{v}))=\left\{b^n,\ b^iab^{n-i-1} :\ i=0,1,\dots, n-1 \right\}$, we obtain $$\chi_n(F(\textbf{v}))=\left\{(0,n), (1,n-1)\right\}.$$ 
Therefore, $\text{c}^{ab}_{F(\textbf{v})}(n)= 2$.

\item If $n_0-l+1\leq n\leq n_0$. Then, we have:
 $$L_n(F(\textbf{v}))=\{b^n,\ b^iab^{n-i-1},\ b^jab^{n_0-l-1}ab^{n-n_0+l-1-j} : i=0,\dots, n-1; \ 0\leq j \leq n-n_0+l-1\}.$$ 
It follows that, $\chi_n(F(\textbf{v}))=\left\{(0,n),\ (1,n-1),\ (2,n-2)\right\}.$ Consequently, $\text{c}^{\text{ab}}_{F(\textbf{v})}(n)= 3$.
 
\item Consider $n> n_0$. According to Theorem \ref{cc}, the classical complexity function of $F(\textbf{v})$ is unbounded. Therefore, $F(\textbf{v})$ is non-ultimately periodic. Thus, we have that $\rho^{ab}_{F(\textbf{v})}(n)\geq 2$. Furthermore, according to Lemma \ref{eq}, the word $F(\textbf{v})$ is $2$-balanced. Since $F(\textbf{v})$ is a binary word, Proposition \ref{pab-binaire} tells us that $\text{c}^{ab}_{F(\textbf{v})}(n)\leq 3$. Therefore, $\text{c}^{ab}_{F(\textbf{v})}(n)\in \left\{2,\hspace{0.1cm} 3 \right\}$.
\end{enumerate}
\end{proof}

\acknowledgements
\label{sec:ack}
We are grateful to Professors Th{\'e}odore Tapsoba and Idrissa Kabor{\'e} for their ongoing contributions to combinatorics on words, as well as for their pioneering work in establishing this field of study in Burkina Faso. We would also like to thank the anonymous reviewers who provided valuable feedback.

\nocite{*}
\bibliographystyle{abbrvnat}
% use the following instead if you encounter problems 
%\bibliographystyle{alpha}
\bibliography{sample-dmtcs}
\label{sec:biblio}

\end{document}